# SHAPE OF TERRITORIES IN SOME COMPETING GROWTH MODELS

By Jean-Baptiste Gouéré

*Université d'Orléans*

We study two competing growth models. Each of these models describes the spread of a finite number of infections on a graph. Each infection evolves like an (oriented or unoriented) first passage percolation process except that once a vertex is infected by type $i$ infection, it remains of type $i$ forever. We give results about the shape of the area ultimately infected by the different infections.

## 1. Introduction and statement of the main results.

1.1. *Introduction.* We study two competing growth models. The first one is the model introduced by Häggström and Pemantle [5]. The second one has been introduced by Deijfen, Häggström and Bagley [1]. Each of these models describes the spread of a finite number of infections on a graph. Each infection evolves like an (oriented or unoriented) first passage percolation process except that once a vertex is infected by type $i$ infection, it remains of type $i$ forever. More explicitly: with each edge $(x,y)$ is associated a positive passage time $\tau(x,y)$; if a vertex $x$ gets infected by type $i$ infection at time $t$ then vertex $y$ gets infected by the same infection at time $t + \tau(x,y)$ except if, at that time, vertex $y$ has already been infected.

Let us denote by $S_t$ the infected territory at time $t$ (without distinguishing between the different types). It is well known that, under some assumptions on the passage times $\tau$, there exists a norm $N$ on $\mathbb{R}^d$ such that $S_t/t$ converges almost surely to the unit ball for norm $N$. One can therefore expect that, if type $i$ infection starts at $x_i$ and if the $x_i$'s are far apart from each other, then the territory infected by type $i$ infection at infinite time looks like the following Voronoï cell:

$$V_i = \{z \in \mathbb{R}^d : \forall j \neq i; N(z - x_i) < N(z - x_j)\}.$$

In this paper, we prove an abstract theorem which is a weak form of such a result. We then apply the theorem to the two above-mentioned models. The proof of the abstract result strongly relies on ideas that appeared in papers by Garet and Marchand [4] and by Hoffman [6, 7]. In these papers, the authors were interested, among other things, in knowing whether the different territories could all be infinite with positive probability or not. We are not aware of any earlier result about the form of the infected territories for the models we study here. Nevertheless, a strong form of such a result can be found in a paper by Pimentel [10] in a two-dimensional setting. The model studied by Pimentel is based on a first passage percolation process on the Delaunay graph of a Poisson point process. His method is very different from ours.

The rest of this section is organized as follows. In the next subsection, we fix some notation. In the following two subsections, we precisely describe the two models and state our results. In the last subsection, we state the abstract result and give the ideas of its proof.

1.2. *Some notation.* For all this paper, we fix an integer $d \geq 2$. Let $\|\cdot\|$ denote the Euclidean norm on $\mathbb{R}^d$. We denote by $|\cdot|$ the canonical Lebesgue measure on $\mathbb{R}^d$ and by $m_1$ the canonical Lebesgue measure on $\mathbb{R}$. We consider these Euclidean spaces as equipped with their natural Borel $\sigma$-algebra. If $C$ and $D$ are two measurable subsets of $\mathbb{R}^d$ and if $|D|$ is positive, we denote by $\underline{\mathrm{dens}}(C|D)$ the lower relative density of $C$ w.r.t. $D$ defined by

$$\underline{\mathrm{dens}}(C|D) = \liminf_{R \to \infty} \frac{|C \cap D \cap B_R|}{|D \cap B_R|},$$

where $B_R$ denotes the Euclidean closed ball of $\mathbb{R}^d$ centered at the origin and with radius $R$. We denote by $\overline{\mathrm{dens}}(C|D)$ the upper relative density of $C$ w.r.t. $D$ defined by

$$\overline{\mathrm{dens}}(C|D) = \limsup_{R \to \infty} \frac{|C \cap D \cap B_R|}{|D \cap B_R|}.$$

When $D = \mathbb{R}^d$, these definitions reduce to the usual definitions of lower and upper density that we denote by $\underline{\mathrm{dens}}(C)$ and $\overline{\mathrm{dens}}(C)$.

1.3. *Model of Häggström and Pemantle.* We begin by recalling the formalism of classical first passage percolation. Let us consider the graph $\mathbb{Z}^d$ obtained by taking $\mathbb{Z}^d$ as vertex set and by putting an edge between two vertices if the Euclidean distance between them is 1. Let us consider a family of nonnegative i.i.d. r.v. $\tau(e)$ indexed by the set of edges $\mathcal{E}$ of the graph. We interpret $\tau(e)$ as the time needed to travel along the edge $e$ (the graph is unoriented). If $a$ and $b$ are two vertices of $\mathbb{Z}^d$, we call path from $a$ to $b$ any finite sequence of vertices $r = (a = x_0, \ldots, x_k = b)$ such that, for all $i \in \{0, \ldots, k-1\}$, the vertices $x_i$ and $x_{i+1}$ are linked by an edge. We denote by $\mathcal{C}(a,b)$ the set of such paths. The time needed to travel along a path $r = (x_0, \ldots, x_k)$ is defined by

$$\widetilde{T}(r) = \sum_{i=0}^{k-1} \tau(x_i, x_{i+1}).$$

We defined $\widetilde{T}(a,b)$, the time needed to go from $a$ to $b$, by

$$\widetilde{T}(a,b) = \inf\{\widetilde{T}(r) : r \in \mathcal{C}(a,b)\}.$$

Let $\tau_1, \ldots, \tau_{2d}$ be $2d$ i.i.d. r.v. admitting the same law as the $\tau(e)$'s. We make the following assumptions:

(1) $$E(\min(\tau_1, \ldots, \tau_{2d})) < \infty$$

and

(2) $$P(\tau_1 = 0) < p_c(d),$$

where $p_c(d)$ denotes the critical probability for bond Bernoulli percolation on $\mathbb{Z}^d$. With these assumptions it is well known that (see, e.g., [8, 9]), for all $x, y$ in $\mathbb{Z}^d$, one has

(3) $$E(\widetilde{T}(x,y)) < \infty$$

and that there exists a norm $N$ on $\mathbb{R}^d$ such that the following convergence holds:

(4) $\widetilde{T}(0,x)N(x)^{-1}$ converges to 1 in $L^1$ as $\|x\|$ goes to infinity.

Fix $x_1, \ldots, x_k \in \mathbb{Z}^d$ ($k \geq 2$). In the competing model of Häggström and Pemantle, for each index $i$, the territory infected by type $i$ infection is the random subset of $\mathbb{Z}^d$ defined by

$$\widetilde{D}_i(x_1, \ldots, x_k) = \{z \in \mathbb{Z}^d : \forall j \in \{1, \ldots, k\} \setminus \{i\}, \widetilde{T}(x_i, z) < \widetilde{T}(x_j, z)\}.$$

We wish to compare this set to the (deterministic) strict Voronoï cells defined for each $i$ by

$$V_i(x_1, \ldots, x_k) = \{z \in \mathbb{R}^d : \forall j \in \{1, \ldots, k\} \setminus \{i\}, N(x_i - z) < N(x_j - z)\}.$$

[We define $V_i(x_1, \ldots, x_k)$ in the same way when the $x_i$'s belongs to $\mathbb{R}^d$ and not necessarily to $\mathbb{Z}^d$.] For every $x \in \mathbb{R}^d$, $\psi(x)$ is defined by

(5) $$\psi(x) \in \mathbb{Z}^d \quad \text{and} \quad x \in \psi(x) + [-1/2, 1/2[^d.$$

We prove the following result:

THEOREM 1.1. *Let $(x_1, \ldots, x_k)$ be a family of distinct points in $\mathbb{R}^d$ ($k \geq 2$). Let $I$ be the set of indices $i \in \{1, \ldots, k\}$ such that the following inequality holds:*

$$\underline{\mathrm{dens}}(V_i(x_1, \ldots, x_k)) > 0.$$

*Let $\varepsilon > 0$. There exists $M > 0$ such that, for all real $R \geq M$, the following assertions hold:*

1. *For all indices $i \in I$, the lower relative density of*
   $$\{z \in \mathbb{Z}^d : P(z \in \widetilde{D}_i(\psi(Rx_1), \ldots, \psi(Rx_k))) \geq 1 - \varepsilon\} + [-1/2, 1/2[^d$$
   *with respect to*
   $$V_i(Rx_1, \ldots, Rx_k)$$
   *is greater than $1 - \varepsilon$.*
2. *The probability of the event*
   $$\{\forall i \in I, \overline{\mathrm{dens}}(\widetilde{D}_i(\psi(Rx_1), \ldots, \psi(Rx_k)) + [-1/2, 1/2[^d | V_i(Rx_1, \ldots, Rx_k))$$
   $$\geq 1 - \varepsilon\}$$
   *is greater than $1 - \varepsilon$.*

One says that coexistence occurs if each territory is infinite. One defines
$$\mathrm{Coex}(x_1, \ldots, x_k) = \bigcap_{i \in \{1, \ldots, k\}} \{\widetilde{D}_i(x_1, \ldots, x_k) \text{ is infinite}\}.$$

As a consequence of Theorem 1.1, we can, for example, get the following coexistence result:

COROLLARY 1.1. *Let $x_1, \ldots, x_k$ be $k$ distinct points of $\mathbb{R}^d$ ($k \geq 2$). Assume that $N(x_i) = 1$ for all indices $i$. Assume that, for all distinct indices $i, j$, $[x_i, x_j]$ contains a point $z$ such that $N(z) < 1$. Let $\varepsilon > 0$. There exists $M > 0$ such that, for all real $R \geq M$, the probability of the event $\mathrm{Coex}(\psi(Rx_1), \ldots, \psi(Rx_n))$ is greater than $1 - \varepsilon$.*

The theorem and the corollary are proved in Section 3.

REMARK. This result relates the number of infections that can coexist to the number of faces (possibly infinity) of the unit ball with respect to the norm $N$.

In the case where the dimension $d$ equals 2, Corollary 1.1 has been proved, among other things, by Hoffman in [7]. Previously to this paper by Hoffman, the positivity of the probability of coexistence for two types of infections had been established by Häggström and Pemantle [5]—in the case where the dimension $d$ equals 2 and the passage times $\tau$ are exponential—and then in the general case independently by Garet and Marchand [4] and by Hoffman [6].

1.4. *Model of Deijfen, Häggström and Bagley.* We begin by recalling the definition of the growth model introduced by Deijfen in [2]. Instead of using the original construction of the process, we use the construction given later in [3] (our presentation is different, but the construction is the same). Let $\nu$ be a probability measure on $]0,\infty[$. Let $\chi$ be a Poisson point process on $\mathbb{R}^d \times \mathbb{R}_+ \times ]0, +\infty[$ whose intensity is the product of the Lebesgue measure on $\mathbb{R}^d \times \mathbb{R}_+$ and of $\nu$ on $]0, +\infty[$. Let us consider the complete directed graph $G$ with vertex set $\mathbb{R}^d$. We associate a time $\tau$ with each edge as follows:

1. For all $x \in \mathbb{R}^d$ we let $\tau(x, x) = 0$.
2. For each point $(X, \widetilde{\tau}, R) \in \chi$ ($X, \widetilde{\tau}$ and $R$, resp., belong to $\mathbb{R}^d$, $\mathbb{R}_+$ and $]0, +\infty[$) and for each vertex $y \in B_R \setminus \{X\}$, we let $\tau(X, y) = \widetilde{\tau}$.
3. For all edges $(x, y)$ to which we have not yet assigned any passage time, we let $\tau(x, y) = +\infty$.

If $a$ and $b$ are two vertices of $G$, we call path from $a$ to $b$ any sequence of vertices $r = (a = x_0, \ldots, x_k = b)$. We denote by $\mathcal{C}(a, b)$ the set of such paths. With each path $r = (x_0, \ldots, x_k)$ we associate a time defined by

$$\widetilde{T}(r) = \sum_{i=0}^{k-1} \tau(x_i, x_{i+1}).$$

If $A$ and $C$ are two subsets of $\mathbb{R}^d$, we define $\widetilde{T}(A, C)$, the time needed to cover $C$ starting from $A$, by

$$\widetilde{T}(A, C) = \sup_{c \in C} \inf\{\widetilde{T}(r) : a \in A, r \in \mathcal{C}(a, c)\}.$$

Notice that we do not have (in general) the equality in law of $\widetilde{T}(A, C)$ and $\widetilde{T}(C, A)$. Let us define

$$S_t = \{x \in \mathbb{R}^d : \widetilde{T}(B, x) \leq t\},$$

where $B$ denotes the unit Euclidean ball centered at the origin. [Notice that, by abuse of notation, we write $\widetilde{T}(B, x)$ instead of $\widetilde{T}(B, \{x\})$.] We make the following assumption:

(6) There exists $a > 0$ such that $\int \exp(ar)\nu(dr)$ is finite.

Under this assumption, the authors of [1] proved the a.s. convergence of $S_t t^{-1}$ toward a deterministic Euclidean ball (whose radius is finite and positive).

In [1], the authors also introduced and studied a related competing growth model. In the particular case we are interested in (each infection evolves with the same velocity), this model can be defined as follows. Fix $x_1, \ldots, x_k \in \mathbb{R}^d$ ($k \geq 2$). Assume that the sets $x_i + B$ are disjoint. At time $t \geq 0$, the territory infected by type $i$ infection is

$$\{z \in \mathbb{R}^d : \widetilde{T}(x_i + B, z) \leq t \text{ and } \forall j \in \{1, \ldots, k\} \setminus \{i\}, \widetilde{T}(x_i + B, z) < \widetilde{T}(x_j + B, z)\}.$$

At time $t = \infty$, the territory infected by type $i$ infection is

$$\widetilde{D}_i(x_1, \ldots, x_n) = \{z \in \mathbb{R}^d : \forall j \in \{1, \ldots, k\} \setminus \{i\}, \widetilde{T}(x_i + B, z) < \widetilde{T}(x_j + B, z)\}.$$

In [3], the authors are interested in the case $k = 2$. They prove that the probability of $\widetilde{D}_1(x_1, x_2)$ and $\widetilde{D}_2(x_1, x_2)$ being both unbounded is positive. They also give such coexistence results when the initial territories are not assumed to be balls.

For each index $i \in \{1, \ldots, k\}$, we wish to compare $\widetilde{D}_i$ to the deterministic Voronoï cell:

$$V_i(x_1, \ldots, x_k) = \{z \in \mathbb{R}^d : \forall j \in \{1, \ldots, k\} \setminus \{i\}, \|z - x_i\| < \|z - x_j\|\}.$$

We prove in this paper the following result:

THEOREM 1.2. *Let $(x_1, \ldots, x_k)$ be a family of distinct points in $\mathbb{R}^d$ ($k \geq 2$). Assume that (6) holds. Let $I$ denote the set of indices $i \in \{1, \ldots, k\}$ such that the following inequality holds:*

$$\underline{\mathrm{dens}}(V_i(x_1, \ldots, x_k)) > 0.$$

*Let $\varepsilon > 0$. There exists $M > 0$ such that, for all real $R \geq M$, the following assertions are satisfied:*

1. 
$$\forall i \in I, \quad \underline{\mathrm{dens}}(\{z \in \mathbb{R}^d : P(z \in \widetilde{D}_i(Rx_1, \ldots, Rx_k)) \geq 1 - \varepsilon\}|$$
$$V_i(Rx_1, \ldots, Rx_k)) \geq 1 - \varepsilon.$$

2. *The probability of the event*

$$\{\forall i \in I, \overline{\mathrm{dens}}(\widetilde{D}_i(Rx_1, \ldots, Rx_k)|V_i(Rx_1, \ldots, Rx_k)) \geq 1 - \varepsilon\}$$

*is greater than $1 - \varepsilon$.*

One says that coexistence occurs if each territory is unbounded. One defines

$$\mathrm{Coex}(x_1, \ldots, x_k) = \bigcap_{i \in \{1, \ldots, k\}} \{\widetilde{D}_i(x_1, \ldots, x_k) \text{ is unbounded}\}.$$

As a consequence of Theorem 1.2, we get the following coexistence result:

COROLLARY 1.2. *Let $x_1, \ldots, x_k$ be $k$ distinct points of the unit Euclidean sphere in $\mathbb{R}^d$ ($k \geq 2$). Let $\varepsilon > 0$. There exists $M > 0$ such that, for all real $R \geq M$, the probability of the event $\mathrm{Coex}(Rx_1, \ldots, Rx_n)$ is greater than $1 - \varepsilon$.*

1.5. *An abstract result and the ideas of the proof.* The results of the previous subsections are applications of an abstract result that we now state and that we prove in Section 2. Fix a norm $N$ on $\mathbb{R}^d$. Let us consider a probability space $(\Omega, \mathcal{F}, P)$ and a family of random variables defined on that space:

$$T = (T(x,y))_{x,y \in \mathbb{R}^d}.$$

If $(x_1, \ldots, x_k)$ is a family of points in $\mathbb{R}^d$ and if $\delta$ is a real we introduce, for every index $i$, the random set

$$D_i^\delta(x_1, \ldots, x_k) = \{z \in \mathbb{R}^d : \forall j \in \{1, \ldots, k\} \setminus \{i\}, T(x_i, z) < T(x_j, z) - \delta\}$$

and the deterministic strict Voronoï cell of $x_i$ defined by

$$V_i(x_1, \ldots, x_k) = \{z \in \mathbb{R}^d : \forall j \in \{1, \ldots, k\} \setminus \{i\}, N(z - x_i) < N(z - x_j)\}.$$

THEOREM 1.3. *Assume that the following assertions hold:*

1. *The map from $\Omega \times \mathbb{R}^d \times \mathbb{R}^d$ to $\mathbb{R}$ defined by $(\omega, x, y) \mapsto T(x,y)(\omega)$ is measurable.*
2. *For all $z$ in $\mathbb{R}^d$, the families of random variables $(T(x,y))_{x,y}$ and $(T(x-z, y-z))_{x,y}$ have the same law.*
3. *For all $x, y$ in $\mathbb{R}^d$, one has $T(x,y) \geq 0$.*
4. *The supremum*

(7) $$\Lambda := \sup_{x \in \mathbb{R}^d : \|x\| \leq 1} E(T(0,x))$$

*is finite.*
5. *For all $x, y, z$ in $\mathbb{R}^d$, one has $T(x,z) \leq T(x,y) + T(y,z)$.*
6. *$T(0,x)N(x)^{-1}$ converges to $1$ in $L^1$ as the norm of $x \in \mathbb{R}^d$ goes to infinity.*

*Let $(x_1, \ldots, x_k)$ be a family of distinct points in $\mathbb{R}^d$ ($k \geq 2$). Let $\delta$ be a real. Let $I$ be the set of indices $i \in \{1, \ldots, k\}$ such that the following inequality holds:*

$$\underline{\mathrm{dens}}(V_i(x_1, \ldots, x_k)) > 0.$$

*Let $\varepsilon > 0$. There exists $M > 0$ such that, for all real $R \geq M$ and for all $i \in I$, the following assertion holds:*

$$\underline{\mathrm{dens}}(\{z \in \mathbb{R}^d : P(z \in D_i^\delta(Rx_1, \ldots, Rx_k)) \geq 1 - \varepsilon\} | V_i(Rx_1, \ldots, Rx_k)) \geq 1 - \varepsilon.$$

REMARK. The convergence assumption requires $T$ to be symmetric enough, but we do not require stronger symmetry conditions [such as, e.g., $T(x,y) = T(y,x)$ a.s. for all vectors $x, y$].

*Ideas of the proof.* The proof strongly relies on two ideas that appeared elsewhere and that we now describe. In order to simplify, let us assume that there are only two sources of infection ($k = 2$) and that one of the sources is the origin. Let us denote the other one by $x$.

• The first idea appeared independently (and in two different forms) in papers by Garet and Marchand [4] and by Hoffman [6]. For all integer $n \geq 0$, let us consider the sum

$$\sum_{i=0}^{n-1} E(T(-x, ix) - T(0, ix)).$$

For each integer $i \in \{0, \ldots, n-1\}$, $T(-x, ix) - T(0, ix)$ is bounded above by $T(-x, 0)$ (by the triangle inequality satisfied by $T$). Therefore, if the norm of $x$ is large enough,

(8) $\quad\quad T(-x, ix) - T(0, ix)$ is roughly bounded above

by $N(x)$ with high probability

(by the convergence of $T$). But for each integer $i \in \{0, \ldots, n-1\}$,

$$E(T(-x, ix) - T(0, ix)) = E(T(0, (i+1)x) - T(0, ix))$$

(by the stationarity of $T$). Therefore the sum equals

$$E(T(0, nx)) - E(T(0, 0)).$$

As a consequence:

(9) The sum is equivalent to $nN(x)$ as $n$ goes to infinity

(by the convergence of $T$). Using (8) and (9) (among other things) one can prove that, for most integers $i \in \{0, \ldots, n-1\}$, $T(-x, ix) - T(0, ix)$ is of order $N(x)$ with high probability. Therefore, most of the points $z$ of $\mathbb{N}x$ are such that $T(-x, z) - T(0, z)$ is of order $N(x)$ with high probability. [In particular, many points in $\mathbb{N}x$ then belong to $D_2(-x, 0)$ with high probability.]

• The second idea appeared in another paper by Hoffman [7]. We exploit the idea in a different way than Hoffman did but the key (simple but powerful) is the same. Assume that $y$ is a vector such that

(10) $\quad\quad\quad\quad\quad\quad N(y - x) < N(y).$

Notice that, for any $z$, $T(y, z) - T(x, z)$ is bounded above by $T(y, x)$ (by the triangle inequality satisfied by $T$). Therefore

(11) $\quad\quad T(y, z) - T(x, z)$ is roughly bounded above

by $N(y - x)$ with high probability,

provided the norm of $y - x$ is large enough (by the convergence and the stationarity of $T$). But, using the first idea one can prove that most of the points $z$ of $-\mathbb{N}y$ are such that

(12) $\qquad T(y, z) - T(0, z)$ is of order $N(y)$ with high probability

(provided the norm of $y$ is large enough). Writing

$$T(x, z) - T(0, z) = (T(y, z) - T(0, z)) - (T(y, z) - T(x, z))$$

and using (10), (11) and (12) one sees that, for such a $z$, $T(x, z) - T(0, z)$ is positive with high probability. To sum up, we have that, under the assumption (10), most of the points $z$ of $-\mathbb{N}y$ belong to $D_2(x, 0)$ with high probability.

• The full proof is given in Section 2. In Section 2.1, we prove that there are enough vectors $y$ satisfying inequalities such as (10) for our purposes. In Section 2.2, we prove results related to the first idea discussed above. In Section 2.3, we prove results related to the second idea and conclude.

**2. Proof of the abstract result.** We work under the assumptions of Theorem 1.3. In this section we need the following definition. If $(x_1, \ldots, x_k)$ is a family of points in $\mathbb{R}^d$ and if $\delta$ is a real we introduce, for every index $i$, the following subset of $\mathbb{R}^d$:

$$V_i^\delta(x_1, \ldots, x_k) = \{z \in \mathbb{R}^d : \forall j \in \{1, \ldots, k\} \setminus \{i\}, N(z - x_i) < N(z - x_j) - \delta\}.$$

2.1. *Geometrical lemmas.* If $K$ is a subset of $\mathbb{R}^d$, we denote by $H(K)$ the set

$$H(K) = \bigcup_{\lambda \geq 1} \lambda K.$$

LEMMA 2.1. *Let $x \in \mathbb{R}^d \setminus \{0\}$ and let $\delta$ be a real. Then:*

1. *The set $V_1^\delta(0, x)$ is stable by homotheties with center $x$ and ratio greater than $1$.*
2. *The set $V_1^\delta(0, x) - x$ is included in the set $V_1^\delta(0, x)$.*

PROOF. In order to simplify the notation, let us define $V := V_1^\delta(0, x)$. Let us check the first item. Let $z, z' \in \mathbb{R}^d$. Assume that $z$ belongs to $V$ and to $[x, z']$. We then have $N(z') \leq N(z' - z) + N(z)$ and then $N(z') < N(z' - z) + N(z - x) - \delta$ and finally $N(z') < N(z' - x) - \delta$. Therefore $z'$ belongs to $V$. This concludes the proof.

Let us check now the second item. Let $z \in V$. Let $w$ denote the middle of $[0, z]$. By the same arguments as previously, one gets that $w$ belongs to $V$. Therefore, by the first item, $z - x$ belongs to $V$. □

LEMMA 2.2. *Let $x \in \mathbb{R}^d \setminus \{0\}$. Then*

$$\underline{\mathrm{dens}}(V_1(0,x)) > 0.$$

*Let $\varepsilon > 0$. There exists $M_1 > 0$, a real $\delta > 0$ and a compact set $K$ included in the unit Euclidean sphere such that, for all $M_2 \geq M_1$, the following relations hold:*

1. $H(M_2 K) \subset V_1(0,x)$;
2. $-M_2 K \subset V_2^\delta(0,x)$;
3. $\underline{\mathrm{dens}}(H(M_2 K)|V_1(0,x)) \geq 1 - \varepsilon$.

PROOF. Let $x \in \mathbb{R}^d \setminus \{0\}$. In order to simplify the notations, let us write $V$ instead of $V_1(0,x)$. As $V$ is open and contains 0, the first item of Lemma 2.1 enables us to conclude that the lower density of $V$ is positive.

For all positive real $R$, we denote by $A_R$ the set of vectors $u$ of the unit Euclidean sphere such that $x + Ru$ belongs to $V$. The sets $A_R$ are open subsets of the unit sphere (since $V$ is an open subset of $\mathbb{R}^d$). By the first item of Lemma 2.1, the family $(A_R)_R$ is nondecreasing w.r.t. inclusion. Let $A$ denote the union of all $A_R$, $R > 0$. Let $\mu$ be the uniform probability measure on the unit Euclidean sphere.

Let $\varepsilon > 0$. Fix $K$, a compact subset of $A$ such that:

$$\frac{\mu(K)}{\mu(A)} \geq 1 - \varepsilon$$

[as $V$, whose Lebesgue measure is infinite, is included in

$$x + \mathbb{R}_+ A = \{x + ra : r \in \mathbb{R}_+, a \in A\},$$

the real $\mu(A)$ is positive]. Fix a real $M_1 \geq 1$ such that the set $A_{M_1}$ contains the set $K$. [We use the fact that $(A_R)_R$ is a nondecreasing family of open sets whose union contains the compact set $K$.] From the inclusion of $x + M_1 K$ in $V$ and from the compactness of $K$, we deduce the existence of a real $\delta > 0$, that we fix, such that $x + M_1 K$ is included in $V_1^\delta(0,x)$.

Let $M_2 \geq M_1$. Write $\widetilde{K} = M_2 K$. Let $\lambda \geq 1$. From the inclusion of $x + M_1 K$ in $V_1^\delta(0,x)$ and from the first item of Lemma 2.1, we deduce the inclusion of $x + \lambda \widetilde{K}$ in $V_1^\delta(0,x)$. By the second item of Lemma 2.1, $\lambda \widetilde{K}$ is included in $V_1^\delta(0,x)$. Therefore $H(\widetilde{K})$ is included in $V_1^\delta(0,x)$ and then in $V$. The first requirement of the lemma is satisfied. From the inclusion of $x + \widetilde{K}$ in $V_1^\delta(0,x)$, we also deduce, by symmetry, the inclusion of $-\widetilde{K}$ in $V_2^\delta(0,x)$, that is, the second requirement of the lemma. As $K$ is included in $A$ and then in the unit Euclidean sphere, the only remaining thing to be proved is the third requirement of the lemma.

First, let us notice the following inclusions:

$$H(\widetilde{K}) \subset V \subset x + \mathbb{R}_+ A.$$

It is therefore sufficient to prove

$$\underline{\mathrm{dens}}(H(\widetilde{K})|x+\mathbb{R}_+A) \geq 1-\varepsilon.$$

But, for all real $R$ large enough,

$$\frac{|H(\widetilde{K}) \cap (x+\mathbb{R}_+A) \cap B_R|}{|(x+\mathbb{R}_+A) \cap B_R|} = \frac{|H(\widetilde{K})R^{-1} \cap B_1|}{|(xR^{-1}+\mathbb{R}_+A) \cap B_1|} = \frac{|H(\widetilde{K})R^{-1} \cap B_1|}{|\mathbb{R}_+A \cap (B_1 - xR^{-1})|}$$

and therefore converges, as $R$ tends to infinity, toward

$$\frac{|[0,1].K|}{|[0,1].A|} = \frac{\mu(K)}{\mu(A)} \geq 1-\varepsilon$$

($[0,1].K$ is the set $\{rk : r \in [0,1], k \in K\}$ and $[0,1].A$ is the set $\{ra : r \in [0,1], a \in A\}$). This concludes the proof. □

LEMMA 2.3. *Let $(x_1, \ldots, x_k)$ be a family of $k$ distinct vectors of $\mathbb{R}^d$ ($k \geq 2$). Assume*

$$\underline{\mathrm{dens}}(V_1(x_1, \ldots, x_k)) > 0.$$

*Let $\varepsilon > 0$. Then there exists a real $\delta > 0$ and a compact set $K$ included in a Euclidean sphere centered at the origin such that:*

1. *The sets $x_1 - K$ and $\{x_1, \ldots, x_k\}$ are disjoint;*
2. *$x_1 + H(K) \subset V_1(x_1, \ldots, x_k)$;*
3. *$x_1 - K \subset \bigcap_{i \neq 1} V_2^\delta(x_1, x_i)$;*
4. *$\underline{\mathrm{dens}}(x_1 + H(K)|V_1(x_1, \ldots, x_k)) \geq 1-\varepsilon$.*

PROOF. Successively apply Lemma 2.2 with $x = x_2 - x_1, \ldots, x = x_k - x_1$. Let us denote by $M_1(2), \ldots, M_1(k)$, by $\delta(2), \ldots, \delta(k)$ and by $\widetilde{K}(2), \ldots, \widetilde{K}(k)$ the positive real numbers and compact sets given by the lemma.

Let $\delta$ be the smallest of the $\delta(i)$'s and let $\widetilde{K}$ be the intersection of the $\widetilde{K}(i)$'s. Fix $M_1$ a real greater than each of the $M_1(i)$'s and such that $x_1 - M_1\widetilde{K}$ contains none of the $x_i$'s. Finally, let $K = M_1\widetilde{K}$, $K(2) = M_1\widetilde{K}(2), \ldots, K(k) = M_1\widetilde{K}(k)$. Items 1, 2 and 3 of the lemma are satisfied. Let us check that the last one is also satisfied.

Let $i \in \{2, \ldots, k\}$. One has

$$\underline{\mathrm{dens}}(H(K(i))|V_1(0, x_i - x_1)) \geq 1-\varepsilon.$$

As $\underline{\mathrm{dens}}(V_1(0, x_i - x_1)) > 0$, translating by $x_1$ does not change the lower density. Therefore,

$$\underline{\mathrm{dens}}(H(K(i)) + x_1|V_1(x_1, x_i)) \geq 1-\varepsilon.$$

We thus have

$$\overline{\mathrm{dens}}(H(K(i))^c + x_1|V_1(x_1, x_i)) \leq \varepsilon$$

and then

$$\overline{\mathrm{dens}}(H(K(i))^c + x_1 | V_1(x_1, \ldots, x_k)) \leq \varepsilon(\underline{\mathrm{dens}}(V_1(x_1, \ldots, x_k) | V_1(x_1, x_i)))^{-1}$$
$$\leq \varepsilon(\underline{\mathrm{dens}}(V_1(x_1, \ldots, x_k)))^{-1}.$$

So,

$$\overline{\mathrm{dens}}(H(K)^c + x_1 | V_1(x_1, \ldots, x_k)) \leq k\varepsilon(\underline{\mathrm{dens}}(V_1(x_1, \ldots, x_k)))^{-1},$$

and then

$$\underline{\mathrm{dens}}(H(K) + x_1 | V_1(x_1, \ldots, x_k)) \geq 1 - k\varepsilon(\underline{\mathrm{dens}}(V_1(x_1, \ldots, x_k)))^{-1}.$$

This ends the proof. □

LEMMA 2.4. *Let $x_1, \ldots, x_k$ be $k$ vectors of $\mathbb{R}^d$ ($k \geq 2$). Assume the existence of a vector $y$ such that, for all $i \in \{2, \ldots, k\}$, $[y, y + x_i - x_1]$ contains a vector $z$ satisfying $N(z) < N(y)$. Then the following inequality holds:*

$$\underline{\mathrm{dens}}(V_1(x_1, \ldots, x_k)) > 0.$$

PROOF.
We first prove that the set

(13) $$W = \bigcap_{i \in \{2, \ldots, k\}} V_2(x_1, x_i)$$

is nonempty. Let $y$ be as stated in the lemma. Let $i \in \{2, \ldots, k\}$. Let $z_i$ be an element of $[y, y + x_i - x_1]$ such that $N(z_i) < N(y)$. Write $z_i = y + \lambda_i(x_i - x_1)$ with $\lambda_i \in [0, 1]$. As $N(z_i) \neq N(y)$, $\lambda_i$ can not equals 0. This enables us to consider the function $\phi : \mathbb{R}^d \to \mathbb{R}^d$ defined by $\phi(u) = \lambda_i^{-1}(u - y) + x_1$. This function maps $y$ to $x_1$, $z_i$ to $x_i$ and 0 to a point that we call $w_i$. Therefore the point $w_i$ satisfies $N(x_i - w_i) < N(x_1 - w_i)$. By the first item of Lemma 2.1, this property is also satisfied by all the points of the half-line

$$L_i := \{x_1 + \alpha(w_i - x_1), \alpha \geq 1\} = \{x_1 - \alpha y, \alpha \geq \lambda_i^{-1}\}.$$

Let now $w$ be in the intersection of the $L_i$'s. This vector belongs to the set $W$ defined by (13).

We now conclude. Let $i \in \{2, \ldots, k\}$. By symmetry, $x_i + x_1 - W$ is included in $V_1(x_1, x_i)$. By the second item of Lemma 2.1, $2x_1 - W$ is also included in $V_1(x_1, x_i)$. Therefore $2x_1 - W$ is included in $V_1(x_1, \ldots, x_k)$. It remains to check that the lower density of $2x_1 - W$ is positive. By the first item of Lemma 2.1, $W$ is stable by all homotheties with center $x_1$ and ratio greater than 1. As $W$ is moreover open and nonempty, one gets that the lower density of $W - x_1$ is positive. The same property therefore holds for $2x_1 - W$. □

Let us state the following immediate consequence of the previous lemma:

LEMMA 2.5. *Let $x_1, \ldots, x_k$ be $k$ distinct points of $\mathbb{R}^d$ ($k \geq 2$). Assume, for all indices $i$, $N(x_i) = 1$. Assume, for all distinct indices $i, j$, that $[x_i, x_j]$ contains a point $z$ such that $N(z) < 1$. Then, for all indices $i$, one has*

$$\underline{\mathrm{dens}}(V_i(x_1, \ldots, x_k)) > 0.$$

PROOF. One can apply Lemma 2.4 to $(x_1, \ldots, x_k)$ (take $y = x_1$). This gives the result for $i = 1$. The proof follows by several other applications of Lemma 2.4. □

2.2. *Competition between two infections on the line joining infection sources.* Let us recall that $m_1$ denotes the canonical Lebesgue measure on $\mathbb{R}$. Let us recall that we work under the assumptions of Theorem 1.3.

LEMMA 2.6. *Let $\varepsilon > 0$. There exists $M > 0$ such that, for all $x$ in $\mathbb{R}^d$ and all real $\lambda$ satisfying*

(14) $$N(x) \geq M \quad \text{and} \quad \lambda \geq M,$$

*one has*

$$m_1(\{\alpha \in [0, \lambda] : E(T(-x, \alpha x) - T(0, \alpha x)) \geq (1 - \varepsilon) N(x)\}) \geq (1 - \varepsilon) \lambda.$$

PROOF. Let $\varepsilon > 0$. Let $M$ be a real such that the following assertions hold (one uses the convergence of $T$ for the first item):

1. For all $x$ in $\mathbb{R}^d$ such that $N(x) \geq M$, one has

$$(1 - \varepsilon) N(x) \leq E(T(0, x)) \leq (1 + \varepsilon) N(x);$$

2. $M \geq 1$ and $M \geq \Lambda \varepsilon^{-1}$ [$\Lambda$ is defined by (7)].

Let now $x$ in $\mathbb{R}^d$ and $\lambda$ a real. One assume that condition (14) holds. Let us consider the real $I$ defined by

$$I = \int_0^\lambda E(T(-x, \alpha x)) \, d\alpha - \int_0^\lambda E(T(0, \alpha x)) \, d\alpha.$$

(One can check that the integrals are finite with Lemma A.1.) Using the stationarity of $T$, one gets

$$I = \int_1^{\lambda+1} E(T(0, \alpha x)) \, d\alpha - \int_0^\lambda E(T(0, \alpha x)) \, d\alpha$$

$$= \int_\lambda^{\lambda+1} E(T(0, \alpha x)) \, d\alpha - \int_0^1 E(T(0, \alpha x)) \, d\alpha.$$

As $\lambda \geq M \geq 1$ and $N(x) \geq M$ one has, for all $\alpha$ in the interval $[\lambda, \lambda + 1]$, the inequality $N(\alpha x) \geq M$ and then the inequality $E(T(0, \alpha x)) \geq (1-\varepsilon)N(\alpha x)$. One therefore has

$$\int_\lambda^{\lambda+1} E(T(0, \alpha x)) \, d\alpha \geq (1-\varepsilon)\lambda N(x).$$

By Lemma A.1 one gets

(15) $$\int_0^1 E(T(0, \alpha x)) \, d\alpha \leq (\|x\| + 1)\Lambda \leq (cN(x) + 1)\Lambda,$$

where $c$ is a fixed positive real such that the inequality $\|\cdot\| \leq cN$ holds. As $1 \leq M \leq N(x)$, $(cN(x)+1)\Lambda$ is bounded above by $(1+c)N(x)\Lambda$. Moreover, as $\Lambda$ is bounded above by $M\varepsilon$ and then by $\lambda\varepsilon$ one gets

$$\int_0^1 E(T(0, \alpha x)) \, d\alpha \leq (1+c)N(x)\lambda\varepsilon.$$

Therefore, one has

(16) $$I \geq (1 - (2+c)\varepsilon)\lambda N(x).$$

Let

$$A = \{\alpha \in [0, \lambda] : E(T(-x, \alpha x) - T(0, \alpha x)) \leq (1 - \sqrt{\varepsilon})N(x)\}.$$

For all $\alpha$ in the interval $[0, \lambda]$ one gets, using successively the triangle inequality satisfied by $T$ and the stationarity of $T$ and then using $N(x) \geq M$,

$$E(T(-x, \alpha x) - T(0, \alpha x)) \leq E(T(-x, 0)) \leq E(T(0, x)) \leq N(x)(1 + \varepsilon).$$

Using the definition of $I$, one then gets

(17) $$\begin{aligned} I &\leq N(x)(m_1(A)(1 - \sqrt{\varepsilon}) + (\lambda - m_1(A))(1 + \varepsilon)) \\ &\leq N(x)(-m_1(A)(\varepsilon + \sqrt{\varepsilon}) + \lambda(1 + \varepsilon)) \\ &\leq N(x)(-m_1(A)\sqrt{\varepsilon} + \lambda(1 + \varepsilon)). \end{aligned}$$

From (16) and (17) one deduces

$$(1 - (2+c)\varepsilon)\lambda N(x) \leq N(x)(-m_1(A)\sqrt{\varepsilon} + \lambda(1 + \varepsilon)),$$

then

$$m_1(A) \leq (3+c)\lambda\sqrt{\varepsilon}$$

and then

$$m_1(\{\alpha \in [0, \lambda] : E(T(-x, \alpha x) - T(0, \alpha x)) \geq (1 - \sqrt{\varepsilon})N(x)\})$$
$$\geq (1 - (3+c)\sqrt{\varepsilon})\lambda.$$

This concludes the proof. $\square$

LEMMA 2.7. *Let $\varepsilon > 0$. There exists $M > 0$ such that, for all $x, y$ in $\mathbb{R}^d$ satisfying*

(18) $\quad N(x) \geq M \quad and \quad E(T(-x,y) - T(0,y)) \geq (1-\varepsilon)N(x),$

*one has*

$$P(T(-x,y) - T(0,y) \geq (1 - 2\sqrt{\varepsilon})N(x)) \geq 1 - 2\sqrt{\varepsilon}.$$

PROOF. Let $\varepsilon > 0$. Fix a real $M > 0$ such that, for all $x$ in $\mathbb{R}^d$ whose norm $N(x)$ is greater or equal to $M$, the following assertions hold (one uses the convergence of $T$):

1. $P(T(0,x) \geq N(x)(1+\varepsilon)) \leq \varepsilon$.
2. $E|T(0,x)N(x)^{-1} - 1| \leq \varepsilon$.

Let $x$ and $y$ be two vectors in $\mathbb{R}^d$ such that condition (18) hold. Let us define an event $G$ by

$$G = \{T(-x, 0) \leq N(x)(1+\varepsilon)\}.$$

Using (18), the triangle inequality and the stationarity satisfied by $T$, one gets

$$\begin{aligned} E((T(-x,y) - T(0,y))1_G) \\ = E(T(-x,y) - T(0,y)) - E((T(-x,y) - T(0,y))1_{G^c}) \\ \geq N(x)(1-\varepsilon) - E(T(-x,0)1_{G^c}) \\ \geq N(x)((1-\varepsilon) - E|T(0,x)N(x)^{-1} - 1| - P(G^c)). \end{aligned}$$

Using Properties 1 and 2 above and using the stationarity of $T$, one therefore gets

$$E((T(-x,y) - T(0,y))1_G) \geq N(x)(1 - 3\varepsilon).$$

The random variable

$$X = N(x)(1+\varepsilon) - (T(-x,y) - T(0,y))1_G$$

is therefore nonnegative (by definition of $G$ and by the triangle inequality satisfied by $T$) and its expectancy is less or equal to $4\varepsilon N(x)$. Therefore, one gets

$$P(X \leq 2\sqrt{\varepsilon}N(x)) \geq 1 - 2\sqrt{\varepsilon},$$

and then

$$P(T(-x,y) - T(0,y) \geq N(x)(1 - 2\sqrt{\varepsilon})) \geq 1 - 2\sqrt{\varepsilon}$$

which concludes the proof. $\square$

LEMMA 2.8. *Let $\varepsilon > 0$. There exists $M > 0$ such that, for all $x$ in $\mathbb{R}^d$ and all real $\lambda$ satisfying*

(19) $$N(x) \geq M \quad \text{and} \quad \lambda \geq M,$$

*one has*

$$m_1(\{\alpha \in [0, \lambda] : P(T(-x, \alpha x) - T(0, \alpha x) \geq (1-\varepsilon)N(x)) \geq 1 - \varepsilon\}) \geq (1-\varepsilon)\lambda.$$

PROOF. This is a consequence of Lemmas 2.6 and 2.7. □

2.3. *Proof of Theorem* 1.3. Let us recall that we work under the assumptions of Theorem 1.3. Theorem 1.3 is a straightforward consequence of the following lemma:

LEMMA 2.9. *Let $(x_1, \ldots, x_k)$ be a family of distinct vectors in $\mathbb{R}^d$ ($k \geq 2$). One assumes*

$$\underline{\mathrm{dens}}(V_1(x_1, \ldots, x_k)) > 0.$$

*Let $\varepsilon > 0$. Let $\delta$ be a real. Then there exists $M > 0$ such that, for all real $R \geq M$, the following property holds:*

$$\underline{\mathrm{dens}}(\{z \in \mathbb{R}^d : P(z \in D_1^\delta(Rx_1, \ldots, Rx_k)) \geq 1 - \varepsilon\} | V_1(Rx_1, \ldots, Rx_k)) \geq 1 - \varepsilon.$$

PROOF. Let us begin by defining $M$. Let us denote by $\alpha$ and $K$ the positive real and the compact set given by Lemma 2.3. One has in particular, for all $x \in K$ and for all $i \in \{2, \ldots, k\}$, the inequality

$$N(x_1 - x - x_i) < N(x_1 - x - x_1) - \alpha.$$

Using the compactness of $K$ one deduces the existence of a real $\eta > 0$, that we fix, such that, for all $x \in K$ and all $i \in \{2, \ldots, k\}$, the following inequality holds:

$$N(x_1 - x - x_i)(1 + \eta) < N(x_1 - x - x_1)(1 - \eta) - \alpha.$$

One deduces the existence of a real $R_2$ such that, for all real $R \geq R_2$, all $x \in K$ and all $i \in \{2, \ldots, k\}$, the inequality

(20) $$N(Rx_1 - Rx - Rx_i)(1 + \eta) < N(Rx_1 - Rx - Rx_1)(1 - \eta) - \delta$$

holds.

By Lemma 2.8 we get a real $R_0$ such that, for all vector $z$ in $\mathbb{R}^d$ and all real $\lambda$ satisfying

$$N(z) \geq R_0 \quad \text{and} \quad \lambda \geq R_0,$$

one has

(21)
$$m_1(\{\alpha \in [0,\lambda]: P(T(-z,\alpha z) - T(0,\alpha z) \geq (1-\eta)N(z)) \geq 1-\varepsilon\})$$
$$\geq (1-\varepsilon)\lambda.$$

Using the convergence of $T$ one gets a real $R_1$ such that, for all vector $z$ in $\mathbb{R}^d$ satisfying $N(z) \geq R_1$, the following inequality holds:

(22)
$$P(T(0,z) \geq N(z)(1+\eta)) \leq \varepsilon.$$

Moreover, as the sets $x_1 - K$ and $\{x_1, \ldots, x_k\}$ are disjoint, there exists a real $R_3$ such that, for all real $R \geq R_3$, all $x \in K$ and all $i \in \{1, \ldots, k\}$, the inequality

(23)
$$N(Rx_1 - Rx - Rx_i) \geq \max(R_0, R_1)$$

holds.

At last, let define $M$ by $M = \max(R_0, R_1, R_2, R_3)$.

We now check that $M$ satisfies the desired property. Fix $R \geq M$. Let $x \in K$, $i \in \{2,\ldots,k\}$ and $\lambda \geq M$. Let us write, for all nonnegative real $\alpha$,

$$T(Rx_i, Rx_1 + \alpha Rx) - T(Rx_1, Rx_1 + \alpha Rx) = A(\alpha, x, i) + B(\alpha, x, i),$$

where

$$A(\alpha, x, i) = T(Rx_1 - Rx, Rx_1 + \alpha Rx) - T(Rx_1, Rx_1 + \alpha Rx)$$

and

$$B(\alpha, x, i) = T(Rx_i, Rx_1 + \alpha Rx) - T(Rx_1 - Rx, Rx_1 + \alpha Rx).$$

By (23), one can use (21) for $z = Rx$. The stationarity of $T$ therefore enables us to get

$$m_1\{\alpha \in [0,\lambda]: P(A(\alpha, x, i) \geq (1-\eta)N(Rx)) \geq 1-\varepsilon\} \geq (1-\varepsilon)\lambda.$$

By (20) one gets

$$m_1\{\alpha \in [0,\lambda]: P(A(\alpha, x, i) > (1+\eta)N(Rx_1 - Rx - Rx_i) + \delta) \geq 1-\varepsilon\}$$
$$\geq (1-\varepsilon)\lambda.$$

By (22) for $z = -Rx_1 + Rx + Rx_i$ [which can be used thanks to (23)], by stationarity of $T$ and by the triangle inequality satisfied by $T$, one gets, for all $\alpha \geq 0$,

$$P(B(\alpha, x, i) \geq -(1+\eta)N(Rx_1 - Rx - Rx_i)) \geq 1-\varepsilon.$$

From the latest two relations one deduces

$$m_1\{\alpha \in [0,\lambda]: P(T(Rx_i, Rx_1 + \alpha Rx) - T(Rx_1, Rx_1 + \alpha Rx) > \delta) \geq 1-2\varepsilon\}$$
$$\geq (1-\varepsilon)\lambda.$$

As this is true for all $i \in \{2, \ldots, k\}$, one gets

$$m_1\{\alpha \in [0, \lambda] : P(Rx_1 + \alpha Rx \in D_1^\delta(Rx_1, \ldots, Rx_k)) \geq 1 - 2k\varepsilon\} \geq (1 - k\varepsilon)\lambda.$$

Let us define the subset $G$ of $\mathbb{R}^d$ by

$$G = \{z \in \mathbb{R}^d : P(Rx_1 + z \in D_1^\delta(Rx_1, \ldots, Rx_k)) \geq 1 - 2k\varepsilon\}.$$

Set $\widehat{x} = x\|x\|^{-1}$ and let $S$ denote the (common) norm of the vectors of $K$. We have proved

$$m_1\{\beta \in [0, RS\lambda] : \beta\widehat{x} \in G\} \geq (1 - k\varepsilon)RS\lambda.$$

We therefore have (when $\varepsilon \leq k^{-1}$)

$$\int_0^{RS\lambda} \beta^{d-1} 1_G(\beta\widehat{x}) \, d\beta \geq \int_0^{m_1\{\beta \in [0, RS\lambda] : \beta\widehat{x} \in G\}} \beta^{d-1} \, d\beta$$

$$\geq \int_0^{(1-k\varepsilon)RS\lambda} \beta^{d-1} \, d\beta$$

$$= (1 - k\varepsilon)^d \int_0^{RS\lambda} \beta^{d-1} \, d\beta.$$

As the previous result holds for all $x$ in $K$, one gets (integrating over $x$ in $K$, which is included in a Euclidean sphere centered at the origin, in a natural way)

(24) $$|B_{RS\lambda} \cap G \cap \mathbb{R}_+ K| \geq (1 - k\varepsilon)^d |B_{RS\lambda} \cap \mathbb{R}_+ K|.$$

As (see Lemma 2.3)

(25) $$\underline{\operatorname{dens}}(x_1 + H(K) | V_1(x_1, \ldots, x_k)) \geq 1 - \varepsilon,$$

$|\mathbb{R}_+ K|$ is positive. Therefore $\underline{\operatorname{dens}}(H(R.K))$ is also positive. In particular $|H(R.K)|$ is positive and $\underline{\operatorname{dens}}(G|H(R.K))$ makes sense. As (24) holds for all $\lambda \geq M$, one deduces

$$\underline{\operatorname{dens}}(G|H(R.K)) \geq (1 - k\varepsilon)^d,$$

that is,

$$\underline{\operatorname{dens}}(\{z \in \mathbb{R}^d : P(z \in D_1^\delta(Rx_1, \ldots, Rx_k)) \geq 1 - 2k\varepsilon\} - Rx_1 | H(R.K))$$
$$\geq (1 - k\varepsilon)^d.$$

Therefore

$$\underline{\operatorname{dens}}(\{z \in \mathbb{R}^d : P(z \in D_1^\delta(Rx_1, \ldots, Rx_k)) \geq 1 - 2k\varepsilon\} | Rx_1 + H(R.K))$$
$$\geq (1 - k\varepsilon)^d.$$

By (25) we get

$$\underline{\operatorname{dens}}(Rx_1 + H(R.K) | V_1(Rx_1, \ldots, Rx_k)) \geq 1 - \varepsilon.$$

Notice that we also have the inclusion of $Rx_1 + H(R.K)$ in $V_1(Rx_1, \ldots, Rx_k)$ (see Lemma 2.3). With new obvious notation, the two previous inequalities and the previous inclusion can be written as follows:

$$\underline{\text{dens}}(X|Y) \geq (1 - k\varepsilon)^d, \qquad \underline{\text{dens}}(Y|Z) \geq 1 - \varepsilon \quad \text{and} \quad Y \subset Z.$$

Using the inclusion, one can write, for all real $a$ large enough,

$$\frac{|X \cap Z \cap B_a|}{|Z \cap B_a|} \geq \frac{|X \cap Y \cap B_a|}{|Y \cap B_a|} \cdot \frac{|Y \cap Z \cap B_a|}{|Z \cap B_a|}.$$

One deduces

$$\underline{\text{dens}}(\{z \in \mathbb{R}^d : P(z \in D_1^\delta(Rx_1, \ldots, Rx_k)) \geq 1 - 2k\varepsilon\}|V_1(Rx_1, \ldots, Rx_k))$$
$$\geq (1 - k\varepsilon)^d (1 - \varepsilon).$$

This concludes the proof of the lemma. □

**3. Proof of Theorem 1.1 and Corollary 1.1.** We use notation and definitions of Section 1.3. We begin by giving precise references and short proofs for (3) and for the following very weak version of (4): for all $x \in \mathbb{Z}^d$,

(26) $\qquad \widetilde{T}(0, kx)k^{-1}$ converges to a finite constant in $L^1$

as the positive integer $k$ goes to infinity.

SKETCH OF THE PROOF OF (3). The proof is sketched in [8], page 135. The idea is the following. One can find $2d$ disjoint paths from 0 to $e_1 = (1, 0, \ldots, 0)$. Using (1) one can then deduce that $E(\widetilde{T}(0, e_1))$ is finite. One can then conclude by subadditivity and symmetry arguments. □

SKETCH OF THE PROOF OF (26). Let $x \in \mathbb{Z}^d$. For all integers $m, n$ we define $X_{m,n}$ by

$$X_{m,n} = \widetilde{T}(mx, nx).$$

We get the desired result by applying Kingman's theorem (we state it in the Appendix) to this family. [Notice that the third condition of Kingman's theorem is satisfied thanks to (3).] □

PROOF OF THEOREM 1.1. We will apply Theorem 1.3 to the family of random variables $T$ defined as follows. Let $U$ be a random variable uniformly distributed on $[-1/2, 1/2[^d$. We assume that $U$ is independent of the random times $\tau(e)$, $e \in \mathcal{E}$. For all $x, y \in \mathbb{R}^d$, one defines $T(x, y)$ by

$$T(x, y) = \widetilde{T}(\widetilde{x}, \widetilde{y}),$$

where, for all vectors $z \in \mathbb{R}^d$, $\tilde{z}$ denotes the unique element of the singleton

$$\mathbb{Z}^d \cap (z - U + [-1/2, 1/2[^d).$$

Let us notice that, if $z$ belongs to $\mathbb{Z}^d$, then $\tilde{z} = z$. Let us check that $T$ satisfies the assumptions of Theorem 1.3.

1. The proof of the measurability is standard.
2. The idea of the proof of the stationarity is the following: the graph $U + \mathbb{Z}^d$ is invariant under the action of the translations of $\mathbb{R}^d$ and $T$ is a factor of the graph. We now give a more detailed proof. Recall that $\mathcal{E}$ denotes the set of edges of $\mathbb{Z}^d$. Let $z \in \mathbb{R}^d$. We define a map $S_z$ from

   $$\Omega = [-1/2, 1/2[^d \times \mathbb{R}_+^{\mathcal{E}}$$

   to itself by

   $$S_z(u, (x_e)) = (u - z - \psi(u - z), (x_{e-\psi(u-z)})),$$

   where $\psi$ is defined by (5). The stationarity of $T$ is a consequence of the following two facts:

   (a) For all $z$ in $\mathbb{R}^d$, $(U, (\tau_e))$ and $S_z(U, (\tau_e))$ have the same law.
   (b) For all $x, y, z \in \mathbb{R}^d$, one has $T(x+z, y+z, U, \tau_e) = T(x, y, S_z(U, \tau_e))$.

3. The nonnegativity is obvious.
4. The finiteness of $\Lambda$ is a consequence of (3).
5. The triangle inequality is satisfied by $\widetilde{T}$ and therefore by $T$.
6. This is a consequence of (4). Nevertheless, as we have found no statement of the well-known result (4), we will give a full proof. [In the relevant literature, the authors usually give proofs of the a.s. convergence instead of the $L^1$ convergence stated in (4). This requires more arguments. By studying these proofs, it is therefore easy to give a proof of (4).] We find it more convenient to directly give a proof of the convergence of $T$ required by Theorem 1.3. Notice that (4) is a straightforward consequence of that convergence.

   We first apply Lemma A.2. This lemma is just a gathering of arguments which are standard in first passage percolation. Notice that the third condition is satisfied thanks to (26) (and thanks to the relation $\tilde{z} = z$ for all $z \in \mathbb{Z}^d$). Applying the lemma, we get the existence of a seminorm $a$ such that

   (27) $$\frac{T(0, x)}{\|x\|} - a\left(\frac{x}{\|x\|}\right)$$

   converges to 0 in $L^1$ as $\|x\|$ tends to $+\infty$. Let $x \in \mathbb{R}^d$ be such that $a(x) = 0$. We wish to prove $x = 0$. Let $i \in \{1, \ldots, n\}$. Let us consider the symmetry $s : \mathbb{R}^d \to \mathbb{R}^d$ define by $s(x_1, \ldots, x_d) = (x_1, \ldots, x_{i-1}, -x_i, x_{i+1}, \ldots, x_d)$.

Notice that $\widetilde{s(x)}$ has the same law as $s(\widetilde{x})$. Notice also that, for every $y \in \mathbb{Z}^d$, $\widetilde{T}(0, y)$ has the same law as $\widetilde{T}(0, s(y))$. One deduces that $T(0, s(x))$ has the same law as $T(0, x)$. Therefore [by (27)], $a(x) = a(s(x)) = 0$. As $a(0, \ldots, 0, 2x_i, 0, \ldots, 0) \leq a(x) + a(s(x))$, one gets that $a(0, \ldots, 0, 2x_i, 0, \ldots, 0) = 0$ and then that $x_i a(0, \ldots, 0, 1, 0, \ldots, 0) = 0$. But by a result of Kesten (Theorem 1.15 in [9]), (2) ensures $a(0, \ldots, 0, 1, 0, \ldots, 0) \neq 0$. One can therefore conclude that, for each index $i$, we have $x_i = 0$. We have proved that $a$ is a norm. By (27), one then sees the required convergence [and therefore (4)] is satisfied with $N := a$.

We now define the value of the real $\delta$ that appears in the statement of Theorem 1.3. Let $\varepsilon > 0$. As $\Lambda$ is finite, by Lemma A.1 we get

$$\sup_{x \in [-1/2, 1/2[^d} E(T(0, x)) < \infty.$$

Therefore, we can fix a real $\delta$ such that, for all $x$ in $[-1/2, 1/2[^d$, the following holds:

$$P(T(0, x) \geq \delta/4) \leq \varepsilon. \tag{28}$$

With such a definition for $\delta$, we have, for all $x, y \in \mathbb{R}^d$, the following inequality:

$$P(|T(x, y) - T(\psi(x), \psi(y))| \geq \delta/2) \leq 2\varepsilon, \tag{29}$$

where $\psi(x)$ is defined by (5). Indeed, by the triangular inequality fulfilled by $T$ and by symmetry of $T$, one has

$$|T(x, y) - T(\psi(x), \psi(y))| \leq T(\psi(x), x) + T(\psi(y), y).$$

Inequality (29) follows by (28) and by stationarity of $T$.

We now prove the first item of Theorem 1.1. Let $(x_1, \ldots, x_k)$ be a family of distinct points in $\mathbb{R}^d$ ($k \geq 2$). Let $\varepsilon > 0$. Let $M$ be the real given by Theorem 1.3. Let $R \geq M$ and $i \in I$. We have

$$\underline{\text{dens}}(\{y \in \mathbb{R}^d : P(y \in D_i^\delta(Rx_1, \ldots, Rx_k)) \geq 1 - \varepsilon\} | V_i(Rx_1, \ldots, Rx_k)) \geq 1 - \varepsilon$$

($D_i^\delta$ and $D$ are defined w.r.t. $T$). But by (29) one has, for any $y$ in $\mathbb{R}^d$, the following inequality:

$$P(\forall j \in \{1, \ldots, k\} : |T(Rx_j, y) - T(\psi(Rx_j), \psi(y))| \leq \delta/2) \geq 1 - 2k\varepsilon.$$

Therefore, one gets that the set

$$\{y \in \mathbb{R}^d : P(y \in D_i^\delta(Rx_1, \ldots, Rx_k)) \geq 1 - \varepsilon\}$$

is included in the set

$$\begin{aligned}(30) \quad & \{\widetilde{y} \in \mathbb{Z}^d : P(\widetilde{y} \in D_i(\psi(Rx_1), \ldots, \psi(Rx_k))) \geq 1 - (2k+1)\varepsilon\} \\ & + [-1/2, 1/2[^d\end{aligned}$$

(by using the fact that any $y$ belongs to $\psi(y) + [-1/2, 1/2[^d)$. Therefore, the set defined by (30) has a lower relative density w.r.t. $V_i(Rx_1, \ldots, Rx_d)$ greater or equal to $1 - \varepsilon$. But, for all $x, y \in \mathbb{Z}^d$, one has $\widetilde{T}(x, y) = T(x, y)$. Therefore in (30) one can replace $D_i$ by $\widetilde{D}_i$ and the first item is proved.

We now prove the second item of Theorem 1.1. Let $i \in I$. By the first item, fix $M$ such that, for all real $R \geq M$, the following inequality holds:

$$\underline{\mathrm{dens}}(A_i | V_i(Rx_1, \ldots, Rx_k)) \geq 1 - \varepsilon, \tag{31}$$

where

$$A_i = \widetilde{A}_i + [-1/2, 1/2[^d$$

and

$$\widetilde{A}_i = \{\widetilde{y} \in \mathbb{Z}^d : P(\widetilde{y} \in \widetilde{D}_i(\psi(Rx_1), \ldots, \psi(Rx_k))) \geq 1 - \varepsilon\}.$$

Fix $R \geq M$. Let us define two random sets by

$$\widetilde{W}_i = \{\widetilde{y} \in \mathbb{Z}^d : \widetilde{y} \notin \widetilde{D}_i(\psi(Rx_1), \ldots, \psi(Rx_k))\}$$

and

$$W_i = \widetilde{W}_i + [-1/2, 1/2[^d.$$

By definition of $A_i$, for all integer $n > 0$ and all $\widetilde{y} \in \mathbb{Z}^d$, one has

$$E|(\widetilde{y} + [-1/2, 1/2[^d) \cap W_i \cap A_i \cap B_n|$$
$$= P(\widetilde{y} \in \widetilde{W}_i)|(\widetilde{y} + [-1/2, 1/2[^d) \cap A_i \cap B_n|$$
$$\leq \varepsilon |(\widetilde{y} + [-1/2, 1/2[^d) \cap A_i \cap B_n|.$$

Therefore, for all integers $n > 0$, one has

$$E|W_i \cap A_i \cap B_n| \leq \varepsilon |A_i \cap B_n| \leq \varepsilon |B_n|.$$

By Fatou's lemma, one therefore gets

$$E(\underline{\mathrm{dens}}(W_i \cap A_i)) \leq \varepsilon.$$

Therefore,

$$P(\underline{\mathrm{dens}}(W_i \cap A_i) \geq \sqrt{\varepsilon}) \leq \sqrt{\varepsilon}$$

and then the event

$$G_i = \{\underline{\mathrm{dens}}(W_i \cap A_i) \leq \sqrt{\varepsilon}\}$$

satisfies $P(G_i) \geq 1 - \sqrt{\varepsilon}$. But when $G_i$ occurs, one has

$$\underline{\mathrm{dens}}(W_i \cap A_i | V_i(Rx_1, \ldots, Rx_k)) \leq \sqrt{\varepsilon} \, \underline{\mathrm{dens}}(V_i(Rx_1, \ldots, Rx_k))^{-1}.$$

But by (31), one has
$$\overline{\operatorname{dens}}(A_i^c | V_i(Rx_1, \ldots, Rx_k)) \leq \varepsilon.$$

Therefore when $G_i$ occurs, one has [using $W_i \subset (W_i \cap A_i) \cup A_i^c$]
$$\underline{\operatorname{dens}}(W_i | V_i(Rx_1, \ldots, Rx_k)) \leq \varepsilon + \sqrt{\varepsilon} \underline{\operatorname{dens}}(V_i(Rx_1, \ldots, Rx_k))^{-1}$$

and then
$$\overline{\operatorname{dens}}(W_i^c | V_i(Rx_1, \ldots, Rx_k)) \geq 1 - \varepsilon - \sqrt{\varepsilon} \underline{\operatorname{dens}}(V_i(Rx_1, \ldots, Rx_k))^{-1}.$$

As $\underline{\operatorname{dens}}(V_i(Rx_1, \ldots, Rx_k)) = \underline{\operatorname{dens}}(R.V_i(x_1, \ldots, x_k)) = \underline{\operatorname{dens}}(V_i(x_1, \ldots, x_k))$, this concludes the proof. □

PROOF OF COROLLARY 1.1. This a consequence of Theorem 1.1 and Lemma 2.5. □

**4. Proof of Theorem 1.2.** We use the notation and conventions of Section 1.4. In particular, $B$ denotes the unit Euclidean ball and $\chi$ denotes the underlying Poisson point process on $\mathbb{R}^d \times \mathbb{R}_+ \times ]0, +\infty[$. We will apply Theorem 1.3 to the family $T$ defined as follows. For all $x, y \in \mathbb{R}^d$, we let
$$T(x, y) = \widetilde{T}(x + B, y + B).$$

LEMMA 4.1. *The following properties hold:*

1. *For all $z \in \mathbb{R}^d$, the families $(T(x+z, y+z))_{x,y}$ and $(T(x,y))_{x,y}$ have the same law.*
2. *For all $x \in \mathbb{R}^d \setminus \{0\}$, the sequence $(T(kx, (k+1)x))_{k \in \mathbb{Z}}$ is stationary and ergodic.*
3. *For all $x, y \in \mathbb{R}^d$, the random variables $T(x,y)$ and $T(y,x)$ have the same law.*
4. *For all $x, y \in \mathbb{R}^d$ such that $\|x\| = \|y\|$, $T(0,x)$ and $T(0,y)$ have the same law.*

PROOF. The properties stated in the lemma are consequences of related properties of the underlying point process $\chi$. □

LEMMA 4.2. *If $A, C$ and $D$ are measurable subsets of $\mathbb{R}^d$, then*
$$\widetilde{T}(A, D) \leq \widetilde{T}(A, C) + \widetilde{T}(C, D).$$

PROOF. Assume that the right-hand side of the inequality stated in the lemma is finite (otherwise the result is obvious). Let $d \in D$ and $\varepsilon > 0$. Fix $c \in C$ and $r_2 \in \mathcal{C}(c,d)$ be such that $\widetilde{T}(r_2) \leq \widetilde{T}(C,D) + \varepsilon$. Now, fix $a \in A$ and $r_1 \in \mathcal{C}(a,c)$ such that $\widetilde{T}(r_1) \leq \widetilde{T}(A,C) + \varepsilon$. If we concatenate $r_1$ and $r_2$ we get an element $r \in \mathcal{C}(a,d)$ such that $\widetilde{T}(r) \leq \widetilde{T}(C,D) + \widetilde{T}(A,C) + 2\varepsilon$. The lemma follows. □

LEMMA 4.3. *Let $A$ and $C$ be two measurable subsets of $\mathbb{R}^d$. We assume that the Lebesgue measure of $A$ is positive and that $C$ is bounded. Then $E(\widetilde{T}(A,C))$ is finite.*

PROOF. The proof is standard. One can proceed as follows. Fix $r > 0$ such that the probability $\nu([2r,+\infty[)$ is positive (this is possible because $\nu(]0,+\infty[) = 1$). Notice the following property.

CLAIM. *Let $a$ be in $\mathbb{R}^d$ and $D \subset \mathbb{R}^d$ be measurable. If the Lebesgue measure of $D \cap (a + B_r)$ is positive then*

$$E(\widetilde{T}(D, D \cup (a + B_r)))$$

*is finite.*

One can prove the claim as follows. Let $U$ be the first $t \geq 0$ such that $\chi$ possesses a point in $D \cap (a + B_r) \times [0,t] \times [2r,+\infty[$. Let $(X, U, R)$ be the point. The law of $U$ is an exponential law with parameter $|D \cap (a + B_r)| \cdot \nu([2r,+\infty[) > 0$. Therefore $E(U)$ is finite. As $X$ belongs to $(a + B_r)$ and as $R$ is greater or equal to $2r$, the set $(a + B_r)$ is contained in the set $X + B_R$. Therefore

$$\widetilde{T}(D, D \cup (a + B_r)) \leq \widetilde{T}(D, D \cup (X + B_R)) \leq U.$$

The claim follows.

As the Lebesgue measure of $A$ is positive and as $C$ is bounded, one can build a finite sequence $a_1, \ldots, a_n$ of vectors in $\mathbb{R}^d$ such that:

1. for each index $i$, the intersection of the set

$$A_{i-1} := A \cup (a_1 + B_r) \cup \cdots \cup (a_{i-1} + B_r)$$

and $(a_i + B_r)$ has a positive Lebesgue measure;
2. the set $A_n := A \cup (a_1 + B_r) \cup \cdots \cup (a_n + B_r)$ contains the set $C$.

As $\widetilde{T}(A,C)$ is bounded above by $\widetilde{T}(A, A_1) + \widetilde{T}(A_1, A_2) + \cdots + \widetilde{T}(A_{n-1}, A_n)$, the lemma follows from the claim. □

By Theorem 1.3, we will get results on $T$, that is, on the $\widetilde{T}(x + B, y + B)$'s. But we are interested in the $\widetilde{T}(x + B, y)$'s. The following lemma will enable us to turn properties on the former into properties on the latter.

LEMMA 4.4. *Let $\varepsilon > 0$. There exists a real $\delta$ such that, for all $x, y \in \mathbb{R}^d$, the following holds:*

$$(32) \quad P(\widetilde{T}(x+B, y) \leq \widetilde{T}(x+B, y+B) \leq \widetilde{T}(x+B, y) + \delta) \geq 1 - \varepsilon.$$

PROOF. We first explicit a realization of $\chi$. Let $\phi$ be a Poisson point process on $\mathbb{R}^d \times \mathbb{R}_+$ whose intensity is the canonical Lebesgue measure. Let us fix a (measurable) enumeration of the points of $\phi$:

$$\phi = \{(X_n, T_n), n \in \mathbb{N}\}.$$

Let $(R_n)_n$ be an independent sequence of i.i.d. r.v. with common distribution $\nu$. Then $\{(X_n, T_n, R_n), n \in \mathbb{N}\}$ is a point process which has the same law as $\chi$. Until the end of this proof, we use this realization of $\chi$ in the definition of the variables $\widetilde{T}$.

Let $x, y \in \mathbb{R}^n$. Let us introduce the following random subset of $\mathbb{N}$:

$$A = \{n \in \mathbb{N} : \|X_n - y\| \leq R_n\}.$$

Let $L$ be the random set of finite sequences of natural integers defined as follows. The sequence $(n_1, \ldots, n_k)$ belongs to $L$ if the following conditions holds:

1. $k \geq 1$.
2. $n_1, \ldots, n_{k-1}$ and $n_k$ are pairwise distinct.
3. $n_1, \ldots, n_{k-1}$ belong to $A^c$, $n_k$ belongs to $A$.
4. $X_{n_1} \in x + B, X_{n_2} \in X_{n_1} + B_{R_{n_1}}, \ldots, X_{n_k} \in X_{n_{k-1}} + B_{R_{n_{k-1}}}$.

If $m$ and $n$ are nonnegative integers, we define variants of $L$ in the following way (in each case, we only point out the differences w.r.t. the definition of $L$):

1. $L_m$: we require in addition that $n_k$ equals $m$.
2. $L^n$: we require in addition that $n$ does not belong to $\{n_1, \ldots, n_k\}$.
3. $L_m^n$: we require in addition that $n_k$ equals $m$ and that $n$ does not belong to $\{n_1, \ldots, n_k\}$.
4. $\widehat{L}_m$: we require in addition that $n_k$ equals $m$ and we drop the requirement $n_k \in A$.

We define $S$ (and in a similar way $S_m, S^n, S_m^n$ and $\widehat{S}_m$) by

$$S = \inf\{T_{n_1} + \cdots + T_{n_k} : (n_1, \ldots, n_k) \in L\}.$$

Notice that, thanks to Lemma 4.3, $\widetilde{T}(x+B, y)$ is a.s. finite. We assume in this step that $y$ does not belong to $x + B$. This ensures

$$\widetilde{T}(x+B, y) = S.$$

Let us define $N$ as the smallest integer $n$ such that $S_n \leq \widetilde{T}(x+B, y) + 1$. We now show that $R_N$ stochastically dominates $\nu$. Let $n$ be a natural integer. Notice that, on the event $\{n \in A\}$, one has:

1. $S_n = \widehat{S}_n$.
2. For all natural integer $k \neq n$: $S_k = S_k^n$.
3. $S = \inf\{S_j, j \in \mathbb{N}\} = \widehat{S}$ where $\widehat{S}$ is defined as the minimum of $\widehat{S}_n$ and $\inf\{S_k^n, k \neq n\}$.

As

(33) $$\{N = n\} \cap \{\widetilde{T}(x+B, y) < \infty\} \subset \{n \in A\},$$

one therefore has

$$\{N = n\} \cap F = F \cap G_n \cap H_n,$$

where

$$F = \{\widetilde{T}(x+B, y) < \infty\},$$
$$G_n = \{n \in A\} = \{R_n \geq \|X_n - y\|\}$$

and

$$H_n = \{\widehat{S}_n \leq \widehat{S} + 1\} \cap \bigcap_{k<n} \{S_k^n > \widehat{S} + 1\}.$$

Let us recall that $P(F) = 1$. Notice that $H_n$ is independent of $R_n$. Let us condition with respect to $\phi$ and denote by $Q$ the resulting random probability. Using the independence properties stated at the beginning of the proof, one gets, for all $r \geq 0$, the following a.s. inequalities [when $Q(N = n)$ is not equal to 0]:

$$Q(R_N \geq r | N = n) = \frac{Q(\{R_n \geq r\} \cap G_n \cap H_n)}{Q(G_n \cap H_n)}$$
$$= \frac{Q(R_n \geq \max(r, \|X_n - y\|))Q(H_n)}{Q(R_n \geq \|X_n - y\|)Q(H_n)}$$
$$\geq Q(R_n \geq r).$$

Therefore, one has a.s.

$$Q(R_N \geq r) \geq Q(R_0 \geq r)$$

and then

$$P(R_N \geq r) \geq P(R_0 \geq r),$$

as desired.

Fix $\Sigma$ (independent of $x$ and $y$) a finite subset of $\mathbb{R}^d$ such that:

1. $B$ is included in $\Sigma + B_{1/2}$;
2. $\Sigma$ is included in $B$.

We claim that if $z \in \mathbb{R}^d$ and $a \geq 1$ are such that $z + B_a$ contains 0, then there exists $s \in \Sigma$ such that $z + B_a$ contains $s + B_{1/2}$. Indeed, in such a case, $a^{-1}z + B$ contains 0. Therefore $a^{-1}z$ belongs to $B$ and there exists $s \in \Sigma$ such that $a^{-1}z$ belongs to $s + B_{1/2}$. For such an $s$, one has $s + B_{1/2} \subset a^{-1}z + B \subset z + B_a$.

We now prove that, with high probability, the territory infected around $y$ at time $\widetilde{T}(x + B, y) + 1$ is not too small. Fix $r \in ]0, 1]$ independent of $x$ and $y$ such that

$$\nu([r, +\infty[) \geq 1 - \varepsilon.$$

We will show

$$(34) \quad P(\exists s \in \Sigma : \widetilde{T}(x + B, y + rs + B_{r/2}) \leq \widetilde{T}(x + B, y) + 1) \geq 1 - \varepsilon.$$

Assume first that $y$ belongs to $x + B$. Then 0 belongs to $xr^{-1} - yr^{-1} + B_{r^{-1}}$ and, by the property previously proved about $\Sigma$ (with $a = r^{-1}$), there exists $s \in \Sigma$ such that $xr^{-1} - yr^{-1} + B_{r^{-1}}$ contains $s + B_{1/2}$. For such a $s$ one has

$$y + rs + B_{r/2} \subset x + B$$

and then

$$\widetilde{T}(x + B, y + rs + B_{r/2}) = 0.$$

Inequality (34) is therefore satisfied in this case.

Assume now that $y$ does not belong to $x + B$. This is the assumption of the step in which $N$ was defined. Recall that $N$ and $\widetilde{T}(x + B, y)$ are a.s. finite. We work on the associated almost sure event. By (33) one gets that 0 belongs to $X_N - y + B_{R_N}$. On the event $\{R_N \geq r\}$, there exists (as above) $s \in \Sigma$ such that $y + rs + B_{r/2}$ is included in $X_N + B_{R_N}$. For such an $s$, one has, on the event $\{R_N \geq r\}$,

$$\widetilde{T}(x + B, y + rs + B_{r/2}) \leq \widetilde{T}(x + B, X_N + B_{R_N}) \leq S_N \leq \widetilde{T}(x + B, y) + 1.$$

As

$$P(R_N \geq r) \geq P(R_0 \geq r) \geq 1 - \varepsilon,$$

(34) is proved.

We now conclude the proof. By Lemma 4.3 we can fix a real $\delta$ (independent of $x$ and $y$) such that

$$P(\widetilde{T}(B_{r/2}, B_2) \leq \delta) \geq 1 - \varepsilon.$$

By stationarity, we then have, for all $s \in \Sigma$,

$$P(\widetilde{T}(y + rs + B_{r/2}, y + rs + B_2) \leq \delta) \geq 1 - \varepsilon$$

and then (as $y + B$ is included in $y + rs + B_2$)

$$P(\widetilde{T}(y + rs + B_{r/2}, y + B) \leq \delta) \geq 1 - \varepsilon.$$

Therefore the probability of the event

$$F = \bigcap_{s \in \Sigma} \{\widetilde{T}(y + rs + B_{r/2}, y + B) \leq \delta\}$$

is greater than or equal to $1 - \text{card}(\Sigma)\varepsilon$. By (34) and the triangular inequality (Lemma 4.2), we therefore have

$$P(\widetilde{T}(x + B, y + B) \leq \widetilde{T}(x + B, y) + 1 + \delta) \geq 1 - (1 + \text{card}(\Sigma))\varepsilon.$$

As the inequality $\widetilde{T}(x + B, y) \leq \widetilde{T}(x + B, y + B)$ is always fulfilled, the lemma is proved. □

The following result is essentially in [2] and [1] but is not explicitly stated (in these papers, the authors prove the almost sure convergence; this requires more arguments). We therefore state the result and provide a short proof.

THEOREM 4.1. *There exists a constant $\mu > 0$ such that $T(0, x)\|x\|^{-1}$ converges to $\mu$ in $L^1$.*

SKETCH OF THE PROOF. Let $x \in \mathbb{R}^d \setminus \{0\}$. For all integers $m, n$ we define $X_{m,n}$ by

$$X_{m,n} = T(mx, nx).$$

The first condition of Kingman's theorem (we state it in the Appendix) is satisfied thanks to Lemma 4.2. The second and forth ones are satisfied because of Lemma 4.1. The third one is satisfied because of Lemma 4.3. We therefore have in particular the convergence in $L^1$ of $T(0, kx)k^{-1}$ toward a finite constant.

We now apply Lemma A.2 in the Appendix. This lemma is just a gathering of arguments which are standard in first passage percolation. Conditions 1 and 2 of the lemma are satisfied thanks to Lemma 4.1. The third condition is a consequence of what we proved in the beginning of the proof. Condition 4 is a consequence of Lemma 4.3 [if $x$ belongs to $B$ then $x + B$ is a subset of $B_2$ and then $T(0, x) \leq \widetilde{T}(B, B_2)$]. Condition 5 is a consequence of Lemma 4.2. Therefore, there exists a seminorm $a$ such that

$$\frac{T(0, x)}{\|x\|} - a\left(\frac{x}{\|x\|}\right)$$

converges to 0 in $L^1$ as $\|x\|$ tends to $+\infty$. By Item 4 of Lemma 4.1, one gets that $a$ is constant on the unit Euclidean sphere. Let us denote by $\mu$ this constant. The only remaining thing to be proved is the inequality $\mu > 0$. This is Proposition 2.1 in [1]. □

PROOF OF THEOREM 1.2. Let us check that $T$ fulfills the assumptions of Theorem 1.3.

1. The proof of the measurability is standard.
2. The stationarity of $T$ is a consequence of Lemma 4.1.
3. The nonnegativity of $T$ is clear.
4. If $x$ belongs to $B$, on has $x + B \subset B_2$ and then $T(0,x) \leq \widetilde{T}(B, B_2)$. Condition 4 is therefore a consequence of Lemma 4.3.
5. The triangular inequality is a consequence of Lemma 4.2.
6. The convergence condition hold with $N = \mu \|\cdot\|$ by Theorem 4.1.

Let $(x_1, \ldots, x_k)$ be a family of distinct vectors in $\mathbb{R}^d$. Let $\varepsilon > 0$. Let $\delta$ be the real given by Lemma 4.4. Let $M$ be the real given by Theorem 1.3. Let $i$ be in $I$ ($I$ is defined in the statement of Theorem 1.2) and $R \geq M$. We have

$$\underline{\mathrm{dens}}(\{y \in \mathbb{R}^d : P(y \in D_i^\delta(Rx_1, \ldots, Rx_k)) \geq 1 - \varepsilon\}|V_i(Rx_1, \ldots, Rx_k)) \geq 1 - \varepsilon,$$

(the $\widetilde{D}_i^\delta$'s are defined w.r.t. $T$). But by (32), we have, for all $y \in \mathbb{R}^d$, the following inequality:

$$P(F_y) \geq 1 - k\varepsilon,$$

where

$$F_y = \{\forall j \in \{1, \ldots, k\}, \widetilde{T}(Rx_j + B, y) \leq \widetilde{T}(Rx_j + B, y + B)$$
$$\leq \widetilde{T}(Rx_j + B, y) + \delta\}.$$

As, for all $y$,

$$F_y \cap \{y \in D_i^\delta(Rx_1, \ldots, Rx_k)\}$$

is included in

$$\{y \in \widetilde{D}_i(Rx_1, \ldots, Rx_k)\},$$

one has

$$\underline{\mathrm{dens}}(\{y \in \mathbb{R}^d : P(y \in \widetilde{D}_i(Rx_1, \ldots, Rx_k)) \geq 1 - (k+1)\varepsilon\}|V_i(Rx_1, \ldots, Rx_k))$$
$$\geq 1 - \varepsilon.$$

The first item of the theorem is proved.

We now prove the second item of Theorem 1.2. Let $i \in I$. By the first item, fix $M$ such that, for all real $R \geq M$, the following inequality holds:

$$\underline{\text{dens}}(A_i | V_i(Rx_1, \ldots, Rx_k)) \geq 1 - \varepsilon,$$

where

$$A_i = \{y \in \mathbb{R}^d : P(y \in \widetilde{D}_i(Rx_1, \ldots, Rx_k)) \geq 1 - \varepsilon\}.$$

Fix $R \geq M$. Let us define

$$W_i = \{y \in \mathbb{R}^d : y \notin \widetilde{D}_i(Rx_1, \ldots, Rx_k)\}.$$

By definition of $A_i$, for all integer $n > 0$, one has

$$E|W_i \cap A_i \cap B_n| \leq \varepsilon |A_i \cap B_n| \leq \varepsilon |B_n|.$$

We conclude as in the proof of the second item of Theorem 1.1. $\square$

## APPENDIX

We begin by a statement of Kingman's theorem.

THEOREM A.1. *Suppose $(X_{m,n}, 0 \leq m < n)$ ($m$ and $n$ are integer) is a family of random variables satisfying:*

1. *For all integers $l, m, n$ such that $0 \leq l < m < n$, one has $X_{l,n} \leq X_{l,m} + X_{m,n}$.*
2. *The distribution of $(X_{m+k,n+k}, 0 \leq m < n)$ does not depend on the integer $k$.*
3. *$E(X_{0,1}^+) < \infty$ and there exists a real $c$ such that, for all natural integer $n$, one has $E(X_{0,n}) \geq -cn$.*

*Then*

$$\lim_{n \to \infty} E(X_{0,n}) n^{-1} \text{ exists and equals } \gamma = \inf_n E(X_{0,n}) n^{-1},$$

$$X := \lim_{n \to \infty} X_{0,n} n^{-1} \text{ exists a.s. and in } L^1 \quad \text{and}$$

$$E(X) = \gamma.$$

*If, for all $k \geq 1$, the stationary sequence $(X_{nk,(n+1)k}, n \geq 1)$ is ergodic, then $X = \gamma$ a.s.*

Let us fix a norm $N$ on $\mathbb{R}^d$. Let $T = (T(x,y))_{x,y \in \mathbb{R}^d}$ be a family of non-negative random variables. We let

$$\Lambda = \sup_{x \in \mathbb{R}^d : \|x\| \leq 1} E(T(0,x)).$$

The following result is very simple:

LEMMA A.1. *Assume that the following conditions hold:*

1. *For all $z$ in $\mathbb{R}^d$, the families of random variables $(T(x,y))_{x,y}$ and $(T(x-z, y-z))_{x,y}$ have the same law.*
2. *For all $x, y, z$ in $\mathbb{R}^d$, one has $T(x,z) \leq T(x,y) + T(y,z)$.*
3. *For all $x, y$ in $\mathbb{R}^d$, one has $T(x,y) \geq 0$.*

*Then, for all $x, y$ in $\mathbb{R}^d$, one has*

$$E(T(x,y)) \leq (\|y-x\|+1)\Lambda.$$

*Moreover for all $x, y, z$ in $\mathbb{R}^d$, one has*

$$E|T(x,y) - T(x,z)| \leq 2(\|y-z\|+1)\Lambda.$$

PROOF. Let us prove the first item. The vector $x - y$ can be written as the sum of $\|x - y\| + 1$ or less vectors of the unit Euclidean ball. Using the stationarity and the triangle inequality satisfied by $T$, one therefore gets:

$$E(T(x,y)) = E(T(0, y-x)) \leq (\|y-x\|+1)\Lambda.$$

The second item is a consequence of the first one because, by the triangle inequality and the nonnegativity of $T$, one has:

$$|T(x,y) - T(x,z)| \leq \max(T(y,z), T(z,y)) \leq T(y,z) + T(y,z). \qquad \square$$

Let $C$ denotes the set of vectors $x$ in $\mathbb{R}^d$ such that $T(0, kx)k^{-1}$ converges in $L^1$ to a finite constant when the integer $k$ goes to infinity. The following result is standard.

LEMMA A.2. *Assume that the following conditions hold:*

1. *For all $z$ in $\mathbb{R}^d$, the families of random variables $(T(x,y))_{x,y}$ and $(T(x-z, y-z))_{x,y}$ have the same law.*
2. *For all $x, y$ in $\mathbb{R}^d$, one has $E(T(x,y)) = E(T(y,x))$.*
3. *The set*

$$\{x\|x\|^{-1}, x \in C \setminus \{0\}\}$$

*is dense in the unit Euclidean sphere $S$.*
4. *$\Lambda$ is finite [$\Lambda$ is defined by (7)].*
5. *For all $x, y, z$ in $\mathbb{R}^d$, one has $T(x,z) \leq T(x,y) + T(y,z)$.*
6. *For all $x, y$ in $\mathbb{R}^d$, one has $T(x,y) \geq 0$.*

*Then there exists a seminorm $a$ on $\mathbb{R}^d$ such that*

$$\frac{T(0,x)}{\|x\|} - a\left(\frac{x}{\|x\|}\right)$$

*converges to 0 in $L^1$ when $\|x\|$ goes to infinity.*

PROOF. If $x \in \mathbb{R}^d$ belongs to $C$, one denote by $a(x)$ the limit of $T(0, kx)k^{-1}$.
Let $x$ be in $C$. Let us prove the following convergence (in $L^1$, with $\lambda$ in $\mathbb{R}$):

$$\lim_{\lambda \to +\infty} T(0, \lambda x)\lambda^{-1} = a(x). \tag{35}$$

If $\lambda$ is a real greater than 1, one has

$$\frac{T(0, \lambda x)}{\lambda} = \frac{T(0, \lfloor \lambda \rfloor x)}{\lfloor \lambda \rfloor} \frac{\lfloor \lambda \rfloor}{\lambda} + \frac{T(0, \lambda x) - T(0, \lfloor \lambda \rfloor x)}{\lambda},$$

where $\lfloor \lambda \rfloor$ denotes the integer part of $\lambda$. By Lemma A.1, the second term converges to 0. But as $x$ belongs to $C$, the first term, and then the sum, converge to $a(x)$.

As a consequence, $C$ is stable by homothety with center 0 and positive ratio. Therefore $C$ is dense in $\mathbb{R}^d$.

Let us now prove that $a$ can be extended into a continuous map from $\mathbb{R}^d$ to $\mathbb{R}$. If $x$ and $y$ are two vectors of $\mathbb{R}^d$ one has, by Lemma A.1,

$$E|T(0, kx)k^{-1} - T(0, ky)k^{-1}| \leq 2(\|kx - ky\| + 1)k^{-1}\Lambda \leq 2(\|x - y\| + k^{-1})\Lambda$$

and then

$$\limsup_{k \to \infty} E|T(0, kx)k^{-1} - T(0, ky)k^{-1}| \leq 2\|x - y\|\Lambda. \tag{36}$$

One deduces that $a$ is $2\Lambda$-Lipschitz (with respect to the Euclidean norm on $\mathbb{R}^d$). This enables us to extend $a$ by continuity on $\mathbb{R}^d$.

Let us now prove that $C = \mathbb{R}^d$. Let $x$ be in $\mathbb{R}^d$. Let $(x_n)_n$ be a sequence of $C$ which converges to $x$. For all integer $n \geq 0$, applying (36) to $x$ and $x_n$, one gets

$$\limsup_{k \to \infty} E|T(0, kx)k^{-1} - a(x_n)| \leq 2\|x - x_n\|\Lambda$$

and then

$$\limsup_{k \to \infty} E|T(0, kx)k^{-1} - a(x)| \leq 2\|x - x_n\|\Lambda + |a(x_n) - a(x)|.$$

Taking limit with respect to $n$, one deduces the desired result.

Let us show that

$$\frac{T(0, x)}{\|x\|} - a\left(\frac{x}{\|x\|}\right)$$

converges to 0 in $L^1$ when $\|x\|$ goes to infinity. Let $(x_n)_n$ be a sequence of vectors whose sequence of norms converge to infinity. To conclude, it suffices to show that one can extract a subsequence $y_n$ such that

$$\frac{T(0, y_n)}{\|y_n\|} - a\left(\frac{y_n}{\|y_n\|}\right)$$

converges to 0 in $L^1$. Fix $(y_n)_n$ a subsequence of $(x_n)_n$ such that $y_n/\|y_n\|$ converges. Let us denote by $y$ the limit. For all integer $n \geq 0$ one has, by Lemma A.1,

$$E\left|\frac{T(0,y_n)}{\|y_n\|} - \frac{T(0,\|y_n\|y)}{\|y_n\|}\right| \leq \frac{2(\|y_n - \|y_n\|y\| + 1)\Lambda}{\|y_n\|}.$$

This upper bound converges to 0. As

$$\frac{T(0,\|y_n\|y)}{\|y_n\|}$$

converges to $a(y)$ and as $a(y_n\|y_n\|^{-1})$ converges to $a(y)$, the result follows.

To conclude, let us check that $a$ is a seminorm. As all $T(x,y)$ are nonnegative, $a$ is nonnegative. By (35) one gets the relation $a(\lambda x) = \lambda a(x)$ for all real $x$ and all nonnegative real $\lambda$. By (35) and by the symmetry of $T$ (Assumption 2 of the lemma) and by stationarity one gets, for all $x$,

$$a(-x) = \lim_{\lambda \to \infty} E(T(0,-\lambda x)\lambda^{-1}) = \lim_{\lambda \to \infty} E(T(0,\lambda x)\lambda^{-1}) = a(x).$$

This enables us to conclude that $a$ is homogeneous. By (35) and using the stationarity and the triangle inequality satisfied by $T$ one gets, for all $x,y$ in $\mathbb{R}^d$,

$$\begin{aligned}
a(x+y) &= \lim_{\lambda \to \infty} E(T(0,\lambda(x+y))\lambda^{-1}) \\
&\leq \lim_{\lambda \to \infty} E(T(0,\lambda x)\lambda^{-1}) + \lim_{\lambda \to \infty} E(T(0,\lambda y)\lambda^{-1}) \\
&= a(x) + a(y).
\end{aligned}$$

The lemma is proved. $\square$

**Acknowledgment.** I would like to thank Olivier Garet for stimulating discussions and for pointing out reference [7].

UNIVERSITÉ D'ORLÉANS
MAPMO-UMR 6628
B.P. 6759
45067 ORLÉANS CEDEX 2
FRANCE
E-MAIL: [jbgouere@univ-orleans.fr](jbgouere@univ-orleans.fr)